\documentstyle[12pt, leqno]{article}

\begin{document}
\def\s{\subseteq}
\def\n{\noindent}
\def\se{\setminus}
\def\dia{\diamondsuit}
\def\la{\langle}
\def\ra{\rangle}

%--------------------------------------------------------------

\title{Cacti with Extremal PI Index}

\footnotetext{
The first author was supported by Project 11571134, 11371162 of NSFC,
 and by the  Self-determined Research Funds of CCNU   from the colleges basic research and operation of
 MOE., the second author was
partially supported by the Summer Graduate Research Assistantship Program of Graduate School of University of Mississippi,
the third author was partially supported by College of Liberal Arts Summer Research Grant of University of Mississippi. }
\author{Chunxiang Wang$^a$, Shaohui Wang$^b$\footnote{  Corresponding author: Shaohui Wang. Emails: C.  Wang (email: wcxiang@mail.ccnu.edu.cn),  S.  Wang (e-mail: shaohuiwang@yahoo.com),   B.  Wei (e-mail:  bwei@olemiss.edu). } ,   Bing Wei$^b$ \\ 
\small\emph {a. School of Mathematics and Statistics, Central China
Normal University, Wuhan, 430079, PRC}\\\small\emph{b. Department of Mathematics, The University of Mississippi, University, MS 38677, USA}}
\date{Accepted by Transactions on Combinatorics, February 2016}
\maketitle

\begin{abstract} 
The vertex PI index $PI(G)  = \sum_{xy \in E(G)} [n_{xy}(x) + n_{xy}(y)]$ is a distance-based molecular structure descriptor,  where $n_{xy}(x)$ denotes the number of vertices which are closer to the vertex $x$ than to the vertex $y$ and which has been the
 considerable research in computational chemistry dating back to Harold Wiener in 1947.  A connected graph is a cactus if  any two of its cycles have at most one common vertex. In this paper,  
we completely determine the extremal graphs with the
largest and smallest vertex PI indices among all the cacti. As a consequence, we obtain
the sharp bounds with corresponding extremal cacti and extend a known result.

\end{abstract}

\section{Introduction}
Let $G$ be a simple connected graph with vertex set $V (G)$ and edge set $E(G)$.
For $x,y \in V(G)$,  the distance $d(x, y)$  is the number of edges in a shortest path connecting  $x$ and $y$.A vertex  is a pendant vertex if its neighborhood contains exactly one vertex.  An edge of a graph is said to be pendant if one of its vertices is a pendant vertex.  An edge $e \in E(G)$ is a cut edge if the graph deleting $e$ contains two components.

A numerical  representation  that can preserve a structural property of a graph  is mathematically defined as a graphic descriptor or a topological index.  
The Wiener index is the oldest and most thoroughly examined  topological index  used  in chemistry. In 1947,  Harold Wiener\cite{w}  applied   Wiener index to determine physical properties of types of Alkanes known  as Paraffins and defined as
 $$W(G) = \sum_{\{x,y\} \subset V(G)} d(x,y).$$ Compared to Wiener index, Szeged index was given  by Klav$\breve{z}$ar and Gutman\cite{ss} in 1996 as follows:
 $$Sz(G) = \sum _{xy \in E(G)}n_{xy}(x)n_{xy}(y),$$  where $n_{xy}(x)$  is the number of vertices $w \in V(G)$ such that $d(x, w) < d(y, w)$, $n_{xy}(y)$  is the number of vertices $w \in V(G)$ such that $d(x, w) > d(y, w)$ and  $w \neq x, y$.
 Currently,   various work relating Wiener index, Sz index 
and their  chemical meaning  have been already studied, referred to the surveys \cite {sz, 2,3,4}.
Based on the considerable success of  Wiener index and Sz index,
   Khadikar\cite{11}  proposed edge Padmakar-Ivan(PI$_e$) index in 2000, which is used in the field of nano-technology,  as follows:
 $$PI_e(G) = \sum_{e=xy \in E(G)} [n_{ex}(e|G) + n_{ey}(e|G)],$$ where $n_{ex}(e|G)$ denotes the number of  edges which are closer to the vertex $x$ than to the vertex $y$, and 
$n_{ey}(e|G)$ denotes the number of  edges which are closer to the vertex $y$ than to the vertex $x$, respectively. The detailed applications of $PI_e$ indices between chemistry and graph theory are investigated in \cite {21}-\cite{25},\cite{11}-\cite{28}. As this definition does not count edges equidistant from both ends of the edge $e = xy$, Khalifeh et al.\cite{5} continued  to  introduce a new PI index of vertex version below: 
$$PI(G) = PI_v(G) = \sum_{xy \in E(G)} [n_{xy}(x) + n_{xy}(y)],$$
where $n_{xy}(x)$ denotes the number of vertices which are closer to the vertex $x$ than to the vertex $y$.
In addition,  there are nice results  regarding vertex PI index in the study of a computational complexity and the intersection between graph theory and chemistry. In \cite{30}, Das and Gutman obtained a lower bound on the vertex PI index of a connected graph in terms of numbers of vertices, edges, pendent vertices, and clique number.  Hoji et al.\cite{35} provided exact formulas for  the vertex PI indices of   Kronecker product  of a connected graph G and a complete graph. Ili$\acute{c}$ and Milosavljevi$\acute{c}$\cite{31} established basic properties of weighted vertex PI index  and proved some lower and upper bounds.   Pattabiraman and Paulraja\cite{32} presented the expressions for vertex PI indices of the strong product of a  graph and the complete multipartite graph.

The synthetic resins\cite{pl}  of plastic materials  is produced by the composition of phenol and benzene with
formaldehyde in  a base. There are no  common edges between independent benzene rings  in the diphenyl ether and the biphenyl. The relation of these rings can be used to partially predict the strength of heat resistance and flame retardancy.
Based on this property, we explore another typle  of graphs:
A graph is a cactus if it is connected and all of its blocks  are either
edges or cycles, i.e., any two of its cycles have at most one common vertex. Denote the cacti of $n$ vertices and $k$ pendent vertices as $\mathcal{C}$$_{n,k}$ with $n \geq k \geq 0$. Let $\lfloor{x}\rfloor$ be the largest integer which is less than or equal to $x$.
 Up to now, many results were obtained concerning  the cacti  between chemistry and graph theory. In \cite{44},  Li and Yang  determined sharp upper and lower bounds of the cacti in $\mathcal{C}$$_{n,k}$ for special chemical indices of Zagreb indices. Feng and Yu\cite{43}   established the cacti in $\mathcal{C}$$_{n,k}$  with the smallest hyper-Wiener indices, which is a renovated version of Wiener index. 
Wang and Tan\cite{42} characterized the extremal cacti having the largest Wiener and hyper-Wiener indices in  $\mathcal{C}$$_{n,k}$.
Wang and Kang\cite{41} found  the extremal   bounds of another chemical index, Harary index, for the cacti  $\mathcal{C}$$_{n,k}$.  Chen\cite{cs} gave the first three smallest Gutman indices among the cacti.

Motivated by the results of chemical indices and their applications, it is worth noting that it may be much interesting
to characterize the cacti in   $\mathcal{C}$$_{n,k}$  with maximum and minimum vertex PI indices.  The concept of vertex PI index yields the following fact.

 {\bf Fact 1  } \emph { Let $G \in $ $\mathcal{C}$$_{n,k}$ with $n \geq k \geq 0$, then\\
(i) If $G$ is $C_3$, $C_4$ or $C_5$, then $PI(G) = 0$, $8$, $10$.\\
(ii) If $G$ is $C_3$ attaching a pendent edge $e$(say $C_3 \cup e$), then $PI(G) = 4$.
}

 In this paper,  we  determine  graphs with the
largest and smallest vertex PI indices in  $\mathcal{C}$$_{n,k}$, and provide the extremal cacti in Figs $1,2$, which  extends Das and Gutman's result\cite{30} by excluding the number of  edges and cliques for the cacti.  Our main results are as follows.(In Figs 1 and 2, $\circ$ means that the vertex maybe exist.)

\begin{figure}[ht]\label{fig 1}
\setlength{\unitlength}{.06in}
\begin{picture}(40,25)(-30,-8)

\put(-5,0){\circle*{0.5}}  
\put(-5,-10){\circle*{0.5}} 
\put(-5,10){\circle*{0.5}}
\put(-6.5,10){\circle*{0.5}}
\put(-8,10){\circle*{0.5}}
\put(-10,10){\circle*{0.5}}
\put(-12.5,10){\circle*{0.5}}
\put(-15,10){\circle*{0.5}}
\put(-8,10){\circle*{0.5}}
\put(-3,10){\circle*{0.5}}
\put(-2,10){\circle*{0.5}}
\put(-1,10){\circle*{0.5}}
\put(0,10){\circle*{0.5}}
\put(2.5,10){\circle*{0.5}}
\put(5,10){\circle*{0.5}}
\put(-10,-10){\circle*{0.5}}
\put(-3,-10){\circle*{0.5}}
\put(-2,-10){\circle*{0.5}}
\put(-1,-10){\circle*{0.5}}
\put(5,-10){\circle*{0.5}}
\put(-2,-3){\circle{0.75}}
\put(1,-6){\circle{0.75}}

\put(-5,0){\line(0,1){10}}  
\put(-5,0){\line(0,-1){10}} 
\put(-5,0){\line(-1,3){3.3}} 
\put(-5,0){\line(-1,2){5}} 
\put(-5,0){\line(-1,1){10}} 
\put(-5,10){\line(-1,0){3}} 
\put(-10,10){\line(-1,0){5}} 
\put(-5,0){\line(1,2){5}} 
\put(-5,0){\line(1,1){10}}
\put( 0,10){\line( 1,0){5}} 
\put(-5,0){\line(-1,-2){5}} 
\put(-5,0){\line(1,-1){10}}

\put(20,0){\circle*{0.5}}  
\put(20,-10){\circle*{0.5}} 
\put(20,10){\circle*{0.5}}
\put(13.5,10){\circle*{0.5}}
\put(12.5,10){\circle*{0.5}}
\put(11.5,10){\circle*{0.5}}
\put(15,10){\circle*{0.5}}
\put(18.5,10){\circle*{0.5}}
\put(10,10){\circle*{0.5}}
\put(17,10){\circle*{0.5}}
\put(22,10){\circle*{0.5}}
\put(23,10){\circle*{0.5}}
\put(24,10){\circle*{0.5}}
\put(25,10){\circle*{0.5}}
\put(27.5,10){\circle*{0.5}}
\put(30,10){\circle*{0.5}}
\put( 15,-10){\circle*{0.5}}
\put(22,-10){\circle*{0.5}}
\put(23,-10){\circle*{0.5}}
\put(24,-10){\circle*{0.5}}
\put(30,-10){\circle*{0.5}}

\put(20,0){\line(0,1){10}}  
\put(20,0){\line(0,-1){10}} 
\put(20,0){\line(-1,3){3.3}} 
\put(20,0){\line(-1,2){5}} 
\put(20,0){\line(-1,1){10}} 
\put(20,10){\line(-1,0){3}} 
\put( 15,10){\line(-1,0){5}} 
\put(20,0){\line(1,2){5}} 
\put(20,0){\line(1,1){10}}
\put( 25,10){\line( 1,0){5}} 
\put(20,0){\line(-1,-2){5}} 
\put(20,0){\line(1,-1){10}} 
\put(10, -15){$\mbox{Fig. 1 }$}

\put(50,0){\circle*{0.5}} 
\put(50,-10){\circle*{0.5}}
\put(50,10){\circle*{0.5}}
\put(47,10){\circle*{0.5}}
\put(45,10){\circle*{0.5}}
\put(40,10){\circle*{0.5}}
\put(47,10){\circle*{0.5}}
\put(52,10){\circle*{0.5}}
\put(53,10){\circle*{0.5}}
\put(54,10){\circle*{0.5}}
\put(55,10){\circle*{0.5}}
\put(60,10){\circle*{0.5}}
\put(45,-10){\circle*{0.5}}
\put(52,-10){\circle*{0.5}}
\put(53,-10){\circle*{0.5}}
\put(54,-10){\circle*{0.5}}
\put(60,-10){\circle*{0.5}}
\put(55,-5){\circle{0.75}}

\put(50,0){\line(0,1){10}}  
\put(50,0){\line(0,-1){10}} 
\put(50,0){\line(-1,3){3.3}} 
\put(50,0){\line(-1,2){5}} 
\put(50,0){\line(-1,1){10}} 
\put(50,10){\line(-1,0){3}} 
\put(45,10){\line(-1,0){5}} 
\put(50,0){\line(1,2){5}} 
\put(50,0){\line(1,1){10}}
\put( 55,10){\line( 1,0){5}} 
\put(50,0){\line(-1,-2){5}} 
\put(50,0){\line(1,-1){10}} 
\put(48, -15){$\mbox{Fig. 2 }$}

\end{picture}

\end{figure}

\vskip 7mm {\bf Theorem 1  } \emph {Let   $G \in$ $\mathcal{C}$$_{n,k} - \{C_3, C_3 \cup e, C_4,C_5\}$ with $n \geq k \geq  0$, then
$PI(G) \leq (n -1+\lfloor  \frac{n-k-1}{3}\rfloor) (n-2),$
where  the equality holds if and only  if $G$ is a tree for $n \leq k+3$ and otherwise,  one of the following statements holds(See Fig. 1):\\
$(i) \;$  All   cycles have length  $4 $ and there are at most  $ k+2$   cut edges.\\
$(ii)$  All  cycles have length  $4$ except one of length $6$ and there are exact $k$ pendent edges.}

\vskip 2mm {\bf Theorem 2  } \emph {
Let   $G \in$ $\mathcal{C}$$_{n,k} - \{C_3, C_3 \cup e, C_4\}$ with $n \geq k \geq  0$, then
$PI(G) \geq  (n-1)(n-2) - 2\lfloor \frac{n-k-1}{2} \rfloor,$
where  the  equality holds if and only if $G$ is a tree for $n \leq k+2$ and otherwise,  all  cycles have length  $3$ and there are at most $k+ 1$  cut edges(See Fig. 2).
}

\section{Main proofs}
Firstly, we provide some lemmas which are important in the proof of our main results.

{ \bf Lemma 1:} Let $G \in$ $\mathcal{C}$$_{n,k}$ and $e \in E(G)$. Then\\
$(i) \; $ $PI(e) \leq n-2$,  the equality holds if  
 $e$ is a  cut edge   or an edge of an even cycle.\\
$(ii)$ If $e$ is an edge of an  odd cycle $C_o$, then $PI(e)   \leq n-3$. Furthermore,  if $G= C_o$, then $PI(e) = n-3$.\\
$(iii)$ For each odd cycle $C$ of $G$,  $PI(C) = (n-2)(|C| -  1) -2$.

{\bf Proof:} Assume that $e= uv  \in E(G)$. Since $PI(e)$ counts at most $n-2$ vertices, then
 $PI(e) \leq n-2$.
If $e$ is a cut edge, then $G-e$ contains two components 	$G_1$ and $G_2$. Thus, all vertices of $G_1$ are closer to one of $\{u,v\}$, say $u$, and all  vertices of $G_2$ are closer to $v$. Thus, $PI(e) =n_e(u)+n_e(v) = n-2$ if $e$ is a cut edge.
Let $C = v_1v_2...v_av_1$ be a cycle of $G$ and $v_lv_l' \in E(C)$. Since $G$ is a cactus, then $G- E(C)$ contains $a$ components $B_1, B_2, ..., B_a$ such that $v_i \in V(B_i)$. If $a$ is even, then $d(v_l, v_i) \neq d(v_l', v_i)$ for   $1 \leq i \leq a$, and $d(v_l, u_i) \neq d(v_l', u_i)$ with $u_i \in V(B_i)$. We obtain that $PI(e) = n-2$ if $C$ is even. Thus, $(i)$ is true.

For $C=C_o$, $a$ is odd. Then there exists a unique vertex $v_t \in V(C)$ such that   $d(v_l, v_t) = d(v_l', v_t)$, that is, $PI(e) \leq n-3$. When $G= C_o$, we see $PI(e) = n-3$. Thus, $(ii)$ is true.

For $(iii)$, $a$ is odd and $\sum_{i=1}^a |B_i| = n$. Note that if $d(v_l, v_t) = d(v_l', v_t)$ with $v_t \in V(C)$, then $d(v_l, u_t) = d(v_l', u_t)$ with $u_t \in V(B_t)$. Similarly,  if $d(v_l, v_t) \neq d(v_l', v_t')$ with $v_t' \in V(C)$, then $d(v_l, u_t') \neq d(v_l', u_t')$ with $u_t' \in V(B_t)$.  Thus, $PI(v_lv_l') = n-2- |B_t|$ with $t \neq l,l'$.  It induces that
$$\begin{array}{rcl}
PI(C) &=& \sum_{e \in E(C)} PI(e) = \sum_{i = 1}^{a} (n-2- |B_i|) \\
&=& a(n-2) - \sum_{i=1}^a|B_i| \\
&=& |C|(n-2) - n \\
&=& (|C|-1)(n-2) -2
\end{array}$$
 and Lemma 1 is true.
$\hfill\Box$

{ \bf Lemma 2:}  Let $C$ be a cycle of $G$. Define 
{\it Transformation 1:} $G_1 = G -xy$ with $xy \in E(G) - E(C)$ and {\it Transformation 2:} $G_2 = G+ x'y'$, where at least one of $\{x',y'\}$ are in $V(G)-V(C)$.
 If $G_1,G_2 \in$ $\mathcal{C}$$_{n,k}$ and $e \in E(C)$, then $PI(e) = PI_{G_1}(e) = PI_{G_2}(e)$.

{\bf Proof:}
Let $C=v_1v_2...v_av_1$, $v_lv_l' \in E(C)$. Then  $G- E(C)$ contains $a$ components $B_1, B_2, ..., B_a$ such that $v_i \in V(B_i)$. Since $G$ is a cactus, then for $v_i  \in V(C)$, if $d(v_l, v_i) =d(v_l', v_i)$, we obtain $d(v_l, u_i) =d(v_l', u_i)$ with $u_i \in V(B_i)$. Similarly,   if $d(v_l, v_i) \neq d(v_l', v_i)$, we obtain $d(v_l, u_i) \neq d(v_l', u_i)$ with $u_i \in V(B_i)$. 
Note that $G_1$ and $G_2$ contain the same cycle $C$ as $G$, and the components $B_j^i$ of $G_i - C$ with
$v_j \in V(B_j^i)$ has the property that $V(B_j^i) = V(B_j^{i'})$.
 Then for $v_i  \in V(C)$, if $d(v_l, v_i) =d(v_l', v_i)$, then $d_{G_1}(v_l, v_i) =d_{G_1}(v_l', v_i)$ and $d_{G_2}(v_l, v_i) =d_{G_2}(v_l', v_i)$, $d_{G_1}(v_l, u_i) =d_{G_1}(v_l', u_i)$ with $u_i \in V_{G_1}(B_i)$ and $d_{G_2}(v_l, u_i) =d_{G_2}(v_l', u_i)$ with $u_i \in V_{G_1}(B_i)$. Similarly,   if $d(v_l, v_i) \neq d(v_l', v_i)$, then $d_{G_1}(v_l, v_i) \neq d_{G_1}(v_l', v_i)$ and $d_{G_2}(v_l, v_i) \neq d_{G_2}(v_l', v_i)$, $d_{G_1}(v_l, u_i) \neq d_{G_1}(v_l', u_i)$ with $u_i \in V_{G_1}(B_i)$ and $d_{G_2}(v_l, u_i) \neq d_{G_2}(v_l', u_i)$ with $u_i \in V_{G_2}(B_i)$. 
 Thus, $PI(e) = PI_{G_1}(e) = PI_{G_2}(e)$ and Lemma 2 is true.
$\hfill\Box$

{ \bf Lemma 3:} If $G \in $ $\mathcal{C}$$_{n,k}$ contains $t_1$ cycles of lengths $\{l_1,l_2,...,l_{t_1}\}$ and $t_2 \geq k$ cut edges, then $PI(G)$ is unique and 
these cycles can be shared a common vertex $u_0$, $k-1$ pendent edges can be adjacent to $u_0$ and a path of length $t_2 - k+1$ can be adjacent to $u_0$. (See Fig. 2) 

{\bf Proof:} By Lemma $1(i)$ and $(iii)$,  PI values with cycles of fixed lengths  and fixd number of cut edges are  determined.  Then $PI(G) = \sum_{C \mbox{ \small is  a cycle of G}} \sum_{e \in E(C)}PI(e) + \sum_{e \mbox{ \small is an cut edge of G}} PI(e)$ is unique.
By recombining these  cycles and cut edges, $t_1$ cycles  can have a common vertex $u_0$, $k-1$ pendent edges can be adjacent to $u_0$ and a path of length $t_2 - k+1$ can be adjacent to $u_0$. Thus, Lemma 3 is true.
$\hfill\Box$

{ \bf Lemma 4:} Let $G \in$ $\mathcal{C}$$_{n,k} -\{C_3 ,C_3\cup e,C_5\}$, if  $PI(G)$ attains the maximal value, then the length of each cycle, if any, is even.

{\bf Proof:}If $G$ has a cycle, then $n \geq 3$.  Assume that there is an odd cycle $C_{2t+1} = u_1u_2...u_{2t}u_{2t+1}u_1$ with  $ t \geq 1$. 
If all vertices of $C_{2t+1}$ have degree 2, then $G = C_{2t+1}$. Since $G \neq C_3, C_5$, then  $n \geq 7$. By Lemma $1(ii)$,  $PI(e) = n-3$ for $e \in E(C_{2t+1})$ and $PI(C_{2t+1}) = n(n-3)$. By Lemma $1(iii)$, $PI(G) = (n-2)(2t)-2$. 
We build a new graph $G' = (G - \{u_1u_{2t+1}\}) \cup \{u_1u_{2t-2}, u_{2t+1}\}$. 
Then $G'$ contains a cycle $C'_1 = u_{2t-2}u_{2t-1}u_{2t}u_{2t+1}u_{2t-2}$ of length 4 and a cycle $C'_2 = u_1u_2...u_{2t-2}u_1$ of length $2t-2$. By Lemma $1(i)$,  $PI(G') = PI(C'_1)+PI(C'_2) = (n-2)(2t+2)$. Thus,
 $PI(G') > PI(G)$, contradicted that $PI(G)$ is maximal.

Thus, there is a vertex of degree at least 3 in   $C_{2t+1}$. If the vertex of degree 3 is unique, say  $u_1$, then there exists a pendent path $u_1v_1v_2...$. 
Set $G_0 = (G - \{u_1u_2\} ) \cup \{u_2v_1\}$, then $G_0 \in$ $\mathcal{C}$$_{n,k} -\{C_3 ,C_3\cup e,C_5\}$. By Lemma 1, we obtain $PI(G_1) > PI(G)$, a contradiction.
If at least two vertices of $\{u_1,u_2,u_3\}$ has degree at least two, say $u_1,u_2$.  Set $G_1 = G -\{u_1u_2\}$, then $G_1 \in$ $\mathcal{C}$$_{n,k} -\{C_3 ,C_3\cup e,C_5\}$. By Lemma 1, we obtain $PI(G) = PI(C) +k(k+1) =  k(k+3)$ and $PI(G_1) = (k+1)(k+3) > PI(G)$, a contradiction. 
If $t \geq 2$, we construct a new graph $G_2$ such that $G_2 = G - \{u_1u_{2t+1}\} \cup \{u_1u_{2t}\}$ with $d_G(u_{2t+1}) \geq 3$. Then $G_2 \in$ $\mathcal{C}$$_{n,k}$, $C_{2t}$ is an even cycle and $u_{2t}u_{2t+1}$ is a cut edge. By Lemma 1 and 2, 
$$\begin{array}{rcl}
PI(G_2) - PI(G) &=& (PI(u_{2t}u_{2t+1}) + PI(C_{2t})) - PI(C_{2t+1}) \\&=& (n-2)(2t+1) - [(n-2)(2t)-2] \\&>& 0,
\end{array}$$  contradicted that $PI(G)$ is maximal. Therefore, each cycle, if any,  is even and Lemma 4 is true.
$\hfill\Box$

{\bf Lemma 5:} Let $G \in$ $\mathcal{C}$$_{n,k} -\{C_3 ,C_3\cup e,C_5\}$ with $n \geq k+4$, if $PI(G)$ attains the maximal value, then all cycles are length $4$ except at most one of them is 6.

{\bf Proof:} By Lemma 4, all cycles are even.
If there exists an cycle $C=u_1u_2...u_{2t}u_1$ with $t \geq 4$. Set $G_1 = (G - \{u_1u_{2t}\}) \cup \{u_1u_4, u_4u_{2t}\}$. Then $G_1 \in$ $\mathcal{C}$$_{n,k} - \{C_3, C_3 \cup e\}$  and $|E(G_1)| = |E(G)| +1$. Since each edge of $G_1$ is either a cut edge or an edge of an even cycle, then $PI(G_1) > PI(G)$ by Lemma $1(i)$, that is, the length of cycles are at most $6$.  
Now suppose that there are two cycles of length 6. By Lemma 3, we can assume these two cycles share a common vertex $u_1$, say $C_1 = u_1u_2...u_6u_1$ and $C_2 = u_1v_2...v_6u_1$. 
 Set $G_2 = G- \{u_1u_2, u_3u_4, u_1v_2\} \cup \{u_1u_4, u_2v_2, u_3v_3,u_1v_3\}$. Then $G_2 \in$ $\mathcal{C}$$_{n,k} - \{C_3, C_3 \cup e\}$  and $|E(G_2)| = |E(G)| +1$. Since each edge of $G_2$ is either a cut edge or an edge of an even cycle, then $PI(G_1) > PI(G)$, that is, there are at most one cycle of length $6$ and Lemma 5 is true.
$\hfill\Box$

{ \bf Lemma 6:} Let $G \in$ $\mathcal{C}$$_{n,k} - \{C_4\}$, if $PI(G)$ attains the minimal value, then the length of each cycle, if any, is odd.

{\bf Proof:} Suppose $G$ has an even cycle $C_{2t} = u_1u_2...u_{2t}u_1$, then  $n \geq k+4$ and $t \geq 2$. 
If all  vertices of $G$ have degree 2, then $G = C_{2t}$ and $n=2t$. By Lemma $1(i)$, $PI(G) = n(n-2) = 2t(2t-2)$. Since $G \neq C_4$ and $t \geq 3$,  set  $G_1 = (G - \{u_{1}u_{2}\}) \cup \{u_1u_{4}, u_2u_{4}\}$.  Then  $G_1 \in $ $\mathcal{C}$$_{n,k} - \{C_4\}$,  $C_{1,3} = u_2u_3u_4u_2$ is an odd cycle and $C_{1, 2t-2} = u_1u_4u_5...u_{2t}u_1$ is an even cycle. By Lemma $1(i)$ and $(iii)$, $PI(G_1) = PI(C_{1, 3}) + PI(C_{1, 2t-2}) = (n-2)2-2 +(n-2)(2t-2) =2t (2t-2) -2 < PI(G),$  contradicted that $PI(G)$ is minimal.
If there exists a vertex $u_2$ with $d(u_2) \geq 3$,  then we construct a new graph $G_2 = (G - \{u_1u_2\}) \cup \{u_1u_3\}$. Then $G_2 \in$ $\mathcal{C}$$_{n,k}$, $u_2u_3$ is a cut edge and $C' = u_1u_3u_4...u_{2t}u_1$ is an odd cycle. By Lemma 1 and 3, 
$$\begin{array}{rcl}
PI(G_2) - PI(G)  &=& (PI_{G_2}(u_2u_3) + PI_{G_2}(C')) - PI(C_{2t})  \\&=&  [(n-2)+(n-2)(2t-2) - 2]- 2t(n-2)  \\&=& -n < 0.
\end{array}$$
Thus, $PI(G_2) < PI(G)$, contradicted that $PI(G)$ is minimal. Therefore, each cycle, if any,  is odd and Lemma 6 is true.
$\hfill\Box$

{ \bf Lemma 7:} Let $G \in$ $\mathcal{C}$$_{n,k} - \{C_4\}$ with $n \geq k+3$, if $PI(G)$ attains the minimal value, then all cycles have length 3. 

{\bf Proof:} By Lemma 6, we only consider all  cycles of $G$ are odd.
Suppose that there is an odd cycle of length greater than 3, say $C_{2t+1} = u_1u_2...u_{2t+1}u_1$ with $t \geq 2$. Set a new graph $G_1 =  (G - \{u_{2t-1}u_{2t}\}) \cup \{u_{1}u_{2t-1}, u_{1}u_{2t}\}$. 
Then $G_1 \in$ $\mathcal{C}$$_{n,k}$ and  we will show that $PI(G_1) < PI(G)$. 
 Let $C_1 = u_1u_2...u_{2t-1}u_1$ and $C_2 = u_1u_{2t}u_{2t+1}u_1$.  By Lemma $1(iii)$, $PI(C) = (n-2)(|C|-2)-2 = 2t(n-2)-2$ and $PI(C_1)+PI(C_2) =[(n-2)(|C_1|-2)-2] + [(n-2)(|C_2|-2)-2] = 2t(n-2)-4$. Thus,  $ PI(C_1)+PI(C_2) < PI(C)$. By Lemma 2, $PI(G_1) - PI(G) = PI(C_1)+PI(C_2) - PI(C) < 0$ and Lemma 7 is true.
$\hfill\Box$

Now, we turn to prove the main results of this paper.

{\bf Proof of Theorem 1.}  
All length of cycles, if any,  are even by Lemma 4. Since $e \in E(G)$ is either a cut edge or an edge of an even cycle, then
  $PI(e) = n-2$ by Lemma $1(i)$.  Thus,  $PI(G) = |E(G)|(n-2)$ and it needs to maximize $|E(G)|$.  For $n \leq k+3$, 
 $\lfloor\frac{n-k-1}{3}\rfloor = 0$ and  $PI(G) = (n-1)(n-2)$.  Thus, Theorem 1 is true.
For $n \geq k+4$,  
all length of cycles are $4$ except at most one of them is 6 by Lemma 5.  By Lemma 3,  all cycles of $G$  have a common vertex $u_0$, $k-1$ pendent edges are adjacent to $u_0$ and a path of length $t_2 - k+1$ is adjacent to $u_0$.

Assume that there exist a cycle $C_6 = u_0u_1u_2u_3u_4u_5u_0$ and $G$ contains more than $k+1$ cut edges, then $G$ has a path $u_0v_1v_2...$ of length more than 2.  Set $G_1 = (G - \{u_2u_3\}) \cup \{u_2v_1, u_0u_3 \}$, then $G_2 \in$  $\mathcal{C}$$_{n,k}$ and $|E(G_2)| = |E(G)| +1$.
 Since  $ e \in E(G_1)$ is either an cut edge or an edge of an even cycle, then $PI(e) = n-2$ and  $PI(G_1) = (n-2)|E(G_1)|> PI(G) = (n-2)|E(G)|$,  contradicted that $PI(G) $ is maximal.  Thus, $G$ contains exact $k$ pendent edges.
 Next we will show that
if all length of cycles are 4, then $G$ contains at most $k+2$ cut edges. Otherwise, there exist a path $u_0v_1v_2...$ of length at least $4$ by Lemma 3. Set $G_2 = G \cup \{u_0v_3\}$, then $G_2 \in$  $\mathcal{C}$$_{n,k}$ and $|E(G_2)| = |E(G)| +1$.  Since  $ e \in E(G_1)$ is either an cut edge or an edge of an even cycle, then $PI(e) = n-2$ and  $PI(G_2) = (n-2)|E(G_2)|> PI(G) = (n-2)|E(G)|$,  contradicted that $PI(G) $ is maximal. 
Note that for $n \geq k+4$,   the number of cycles of $G$ is $\lfloor\frac{n-k-1}{3}\rfloor$ and the number of edges of $G$ is $n -1 + \lfloor\frac{n-k-1}{3}\rfloor$. Thus, $PI(G) = (n -1 + \lfloor\frac{n-k-1}{3}\rfloor)(n-2)  $ and
 Theorem 1 is true.
$\hfill\Box$

{\bf Proof of Theorem 2.} For $n \leq k+2$,  $\lfloor \frac{n-k-1}{2} \rfloor = 0$ and $PI(G) = (n-1)(n-2)$ by Lemma $1$. Thus, Theorem 2 is true.  For $n \geq k+3$, the length of each edge of $G$ is 3 by Lemma 7.  Next we will show that $G$ contains at most $k+1$ cut edges. Assume that $G$ contains at least $k+2$ cut edges. By Lemma 3,  all cycles of $G$  have a common vertex $u_0$, $k-1$ pendent edges are adjacent to $u_0$ and a path of length at least $(k+2) - k+1 = 3$ is adjacent to $u_0$. Denote the path as $u_0v_1v_2v_3....$, set $G_1 = G \cup \{u_0v_2\}$. By Lemma $1(iii)$ and $2$, $PI(G_1) - PI(G) = PI_{G_1}(v_0u_1u_2v_0) - PI(u_0v_1) - PI(v_1v_2) =  [(n-2)(3-1)-2] - (n-2) -(n-2) = -2 <0$. Thus,  
 $PI(G_1) <   PI(G)$, contradicted that $PI(G)$ is minimal.
Note that for $n \geq k+3$, the number of cycles of length 3 is $\lfloor \frac{n-k-1}{2} \rfloor  $ and the number of cut edges is $n-1 - 2\lfloor \frac{n-k-1}{2} \rfloor   $.
Thus, 
$$\begin{array}{rcl}
 PI(G) &=& 2(n-3) (\lfloor \frac{n-k-1}{2} \rfloor) +( n -1- 2\lfloor \frac{n-k-1}{2} \rfloor ) (n-2)\\ 
&=&  (n-1)(n-2) - 2\lfloor \frac{n-k-1}{2} \rfloor,
\end{array}$$
  and Theorem 2 is true.
$\hfill\Box$

{\bf Remarks.} 
The maximal and minimal values of vertex PI vertices of cacti are uniqe, but the cacti achieved the maximal and minimal vertex PI index are not unique. All cacti satisfying the statements in Theorem 1 and Theorem 2 are arrived at the corresoponding sharp values. Fig 1 and Fig 2 are special examples achieved the sharp bounds.

\end{document}